\documentclass[11pt]{article}
\usepackage{amsthm, amssymb, amsmath, geometry, verbatim}
\usepackage{authblk}

\numberwithin{equation}{section}
\usepackage{color}
\usepackage{amsfonts,latexsym,amsmath,amscd,geometry,amssymb}
\newcommand \nc{\newcommand}
\newtheorem{theorem}{Theorem}[section]
\newtheorem{lemma}[theorem]{Lemma}

\newtheorem{proposition}[theorem]{Proposition}
\newtheorem{corollary}[theorem]{Corollary}

\newtheorem{definition}[theorem]{Definition}

\newtheorem{remark}[theorem]{Remark}
\nc{\be}{\begin{eqnarray}}\nc{\ee}{\end{eqnarray}}
\nc{\beq}{\begin{equation}}\nc{\eeq}{\end{equation}}
\nc{\bex}{\begin{eqnarray*}}\nc{\eex}{\end{eqnarray*}}
\nc{\beqno}{\begin{eqnarray*}}\nc{\eeqno}{\end{eqnarray*}}
\nc{\dt}{\frac{\mathrm{d}}{\mathrm{d}t}}

\def\pf{\noindent{\bf Proof.\quad}}\def\endpf{\hfill$\Box$}
\def\di{\mathrm{div\,}}

\allowdisplaybreaks

\usepackage[colorlinks,linkcolor=blue]{hyperref}

\begin{document}
\title{Global well-posedness of 3D   two-fluid type model with vacuum: smallness on scaling invariant quantity}
\author{Huanyao Wen\thanks{Email: mahywen@scut.edu.cn}\quad \quad  Chanxin Xie\thanks{Email: chanxinxie@163.com}}
\affil{School of Mathematics, South China University of Technology, Guangzhou 510641, China.}

\date{}
\maketitle
\begin{abstract}
This paper focuses on Cauchy problem for the three-dimensional 
 two-fluid type model, in which the presence of vacuum is permitted.
Under some assumptions that the initial data satisfy appropriate regularity conditions and a compatibility constraint, and that the newly introduced scaling-invariant initial quantities $\bar P^{\frac{ 3}{\gamma}} \left(\|\sqrt{\rho_0}u_0\|_{L^2}^2+\|P_0\|_{L^1}\right) \left(\|\nabla u_0\|_{L^2}^2+\|P_0\|_{L^2}^2\right)$ and $\bar P^{\frac{6}{\gamma}+1} \left(\|\sqrt{\rho_0}u_0\|_{L^2}^2+\|P_0\|_{L^1}\right)^3 \left(\|\nabla u_0\|_{L^2}^2+\|P_0\|_{L^2}^2\right)$
are sufficiently small, the global well-posedness of strong solutions to the  two-fluid type model  is derived.
\end{abstract}

 \bigbreak \textbf{{\bf Key Words}:} Viscous  two-fluid  model, vacuum, global strong solutions, scaling invariance.

 \bigbreak  {\textbf{AMS Subject Classification 2020:}  76N10, 76T10, 35K65.


 \section{Introduction}
\setcounter{equation}{0}\setcounter{theorem}{0}
\renewcommand{\theequation}{\thesection.\arabic{equation}}
\renewcommand{\thetheorem}{\thesection.\arabic{theorem}}
\subsection{Background and motivation}
\quad \quad Two-fluid flows are ubiquitous in nature and arise in a wide range of industrial applications,
including power generation, nuclear engineering, chemical processing, oil and gas production,
cryogenics, aerospace engineering, biomedical applications, and micro-scale technologies.
They are also relevant to the study of various mathematical models, such as cancer cell migration models \cite{Evje2017,EvjeWen2018},  MHD
system \cite{LiSun2019}, and Oldroyd-B model \cite{BarrettLuSuli2017}. 
 In this paper, we investigate a viscous compressible  two-fluid type model as follows:
\begin{align}\label{twofluid}
\begin{cases}
\begin{aligned}
&\partial_t n + \nabla \cdot (n u)= 0, \\
 &\partial_t \rho + \nabla \cdot (\rho u)= 0, \\
 &\partial_t (\rho u)+ \di  [(\rho u \otimes u)] + \nabla P(n,\rho)= \mu \Delta u + (\mu + \lambda) \nabla \di  u,
\end{aligned}
\end{cases}
\end{align} where $\rho,\,n,\,u$  represent the densities of two
fluids and the velocity of the fluids,\, respectively.  The viscosity coefficients $\mu$ and $\lambda$ represent
viscosities of mixture of two fluids and are assumed to be constants satisfying $\mu>0,\,\, 2\mu+3\lambda \geq 0$ for three-dimensional case.

The mathematical analysis of viscous two-fluid flow models has received considerable attention over the past decades with substantial progress in wellposedness theory. However, there are some fundamental problems still open. In the following, we provide a brief overview of results concerning the model \eqref{twofluid}. Evje and Karlsen \cite{EK2008} first established the global existence theory of weak solutions with arbitrarily large initial data in a one-dimensional regime that excludes transitions to single-phase flow, and hence vacuum is not allowed. The pressure function satisfies
\be\label{pressure1}
P(n,\rho)=C_0\Big[-b(\rho,n)+\sqrt{b(\rho,n)^2+c(\rho,n)}\Big],
\ee where $b(\rho,n)=\kappa_0-\rho-a_0n,\,\,c(\rho,n)=4\kappa_0a_0n$ with positive constants $\kappa_0$ and $a_0$. This result was subsequently extended by Evje, Wen, and Zhu \cite{EvjeWenZhu2017} to allow transition to single-phase flow. See \cite{Yao-Zhu2} for free boundary value problem.
For multi-dimensional case, there are more challenges arising from the two-component pressure. In this context, another pressure law was considered, namely,
\be\label{pressure2}
P(n,\rho)=\rho^\gamma+n^\Gamma
\ee
for $\gamma,\Gamma\ge1$. Vasseur, Wen, and Yu  \cite{VWY-2019} proved the global existence of weak solution in three dimensions for arbitrarily large initial data for $\Gamma > \frac{9}{5}$ and $\gamma \geq 1$ under a density domination condition. Novotn\'y and Pokorn\'y \cite{NP2020} extended this result to $\Gamma=\frac{9}{5}$ and more general pressure functions are considered. Wen \cite{W2021} removed the density domination assumption, allowing both adiabatic exponents to attain the critical value $\frac{9}{5}$. See \cite{KMN2024} for further progress. The uniqueness of the weak solutions is still a challenging open problem. The two-fluid type model (\ref{twofluid}) with pressure (\ref{pressure2}) can be considered as a reduced model from the compressible diffusive Oldroyd-B-type model for the cases $\Gamma = 1,\,2$ \cite{BarrettLuSuli2017} and from compressible isentropic Magnetohydrodynamics(MHD) model \cite{jiangzhang2017,LiSun2019} for the case $\Gamma =2$ in two dimensions where $n$ represents the vertical magnetic field component. See Hu \cite{Hu}, Wu and Wu \cite{Wuwu2017} for further investigation. When the initial data are close to a constant equilibrium, the global wellposedness of strong solutions to (\ref{twofluid}) and (\ref{pressure2}) in three dimensions with small initial energy and a domination condition of densities was proved by Guo, Yang, and Yao \cite{GuoYangYao2011}, and later by Hao and Li \cite{HaoLi2012} in homogeneous Besov spaces. Yu \cite{Yu2021} further extended these results by allowing vacuum regions in both phases and deriving higher-order estimates without the domination condition of densities. However, the smallness of the initial data is not scaling invariant. It is worth noticing that system \eqref{twofluid} is invariant under the following scaling transformation
\begin{equation}\label{scl-inva-non}
\begin{split}
\rho_\lambda(x, t) &= \lambda^{\frac{1}{\gamma}}\rho(\lambda^{\frac{\gamma+1}{2\gamma}} x, \lambda t),\\
  n_\lambda(x, t)& =\lambda^{\frac{1}{\Gamma}} n(\lambda^{\frac{\gamma+1}{2\gamma}} x, \lambda t),\\
u_\lambda(x, t) &= \lambda^{\frac{\gamma-1}{2\gamma}} u(\lambda^{\frac{\gamma+1}{2\gamma}} x, \lambda t).
\end{split}
\end{equation} It means that if $(\rho,n,u)$ is a solution to (\ref{twofluid}), so is $(\rho_\lambda,n_\lambda,u_\lambda)$. (\ref{twofluid}) is complemented by the following initial  condition
\beq\label{non-initial}
\left (\rho,n,u)\right|_{t=0} = ( \rho_0,\,n_0,\,u_0),
\eeq
and the following far-field behavior
\beq\label{non-cauchy}
(\rho,n,u) \to (0,\,0,\,0), \quad \quad \text{as $|x| \to \infty$},
\eeq

 Our main aim is to establish the global well-posedness theory for the three-dimensional compressible two-fluid type system \eqref{twofluid} and (\ref{pressure2}) with the initial conditions (\ref{non-initial}) and (\ref{non-cauchy}) allowing vacuum, under a scaling invariant smallness assumption. In contrast to the single-phase compressible Navier--Stokes equations, the  two-fluid  type system is intrinsically inhomogeneous and strongly coupled through two density variables, which makes the existence of nontrivial scaling invariance far from obvious. In particular, the degeneracy of the equations near vacuum is significantly amplified by the interaction between the liquid and gas phases. Moreover, due to the strong coupling among $(\rho,\, n,\, u)$, the initial data themselves cannot be expected to be scale-invariant, which precludes a direct analysis in critical spaces.

To overcome these difficulties, we identify a suitable scaling-invariant combination of the initial data that reflects the  two-fluid structure and allows uniform control of the nonlinear
terms. More precisely, we introduce the following small scaling-invariant quantity:
\begin{equation}\label{scaling-invariant}
\bar P^{\frac{ 3}{\gamma}} E_0 \big(1+\bar P^{\frac{ 3}{\gamma}}\bar PE_0^2\big)\big(\|\nabla u_0\|_{L^2}^2+\|P_0\|_{L^2}^2\big)
\end{equation}
which plays a fundamental role in closing the global \emph{a priori} estimates.

Motivated by \cite{Wen2025}, we  derive a  uniform bound of
\begin{equation}\label{pressure-estimate}
\int P^2 + \int_0^t \int P^3\, ds,
\end{equation}
which plays a crucial role in controlling the supercritical convective term.
In contrast to the single-phase case, the pressure $P(\rho,n)$ depends
simultaneously on two density variables through the continuity equations.
Such a coupling leads to a significantly more involved nonlinear structure
and generates additional interaction terms in the pressure estimates. Specifically, $P$ satisfies the following equation
\beq
\begin{split}
P_t+u\cdot \nabla P+\gamma P\di u+(\Gamma-\gamma)n^\Gamma\di u=0.
\end{split}
\eeq
We observe that the equation of \( P \) depends not only on $P$ and \( u \), but also on the term \( n^\Gamma \). It turns out that it needs to estimate the term
\[
\int n^\Gamma P \, \di u,
\]
which presents a new challenge. To address this, we employ the effective viscous flux formulation $G=(2\mu+\lambda)\di u-P$ and combine it with the equation for \( n^\Gamma \). On the other hand, we need to derive the $L^\infty$ bounds for $\rho$ and $P$. Given that $G$ is dependent on $P(\rho,n)$, the upper bounds of $\rho$ and $n$ are mutually constrained. Consequently, the estimate $\|\rho\|_{L^\infty}$ cannot be established in the absence of a density domination condition. To resolve this difficulty,  we make a priori assumption for the pressure $P$ instead of the densities, and use the upper bound of pressure to control the densities, see Proposition \ref{prop1}. Specifically, we make use of the relation
\[
\|\rho\|_{L^\infty} \leq \|P\|_{L^\infty}^{\frac{1}{\gamma}},
\]
which enables us to achieve our main results without imposing any density domination condition.

\medskip

The rest of the paper is organized as follows. In Section \ref{sec1.2}, we present main result concerning the global well-posedness of strong solutions to system \eqref{twofluid}. Section \ref{sec 2} is devoted to the proof of the main result.

\bigskip

Before stating our main result, we would like to introduce some notations and definitions which will
be used throughout this paper.
\begin{itemize}
\item $\displaystyle\int  f =\int_{\mathbb{R}^3} f \,{\rm d}x$.
\item  For $1\le l\le \infty$, we use the following notations for the standard Lebesgue and Sobolev spaces:
$$\displaystyle L^l=L^l(\mathbb{R}^3), \ \ \ D^{k,l}=\left\{ u\in L^1_{\rm{loc}}(\mathbb{R}^3): \|\nabla^k u \|_{L^l}<\infty\right\},$$
$$W^{k,l}=L^l\cap D^{k,l}, \ \ \ H^k=W^{k,2}, \ \ \ D^k=D^{k,2},$$
$$D_0^1=\Big\{u\in L^6:  \, \|\nabla u\|_{L^2}<\infty \},$$
$$\|u\|_{D^{k,l}}=\|\nabla^k u\|_{L^l}.$$

\end{itemize}
\subsection{Main result}\label{sec1.2}
We first specify the definition of strong solutions in this paper.
\begin{definition}
Given a time $T>0$, a triple $(\rho,n,u)$ is called a \emph{strong solution}
to the initial-value problem for system \eqref{twofluid} and (\ref{pressure2}) with the initial conditions (\ref{non-initial}) and (\ref{non-cauchy}) in $\mathbb{R}^3\times[0,T]$ if
\begin{align*}
&(\rho,n)\in C\big([0,T];W^{1,q}(\mathbb{R}^3)\cap H^1(\mathbb{R}^3)\big),\quad
\rho,n\ge0,\\
&(\rho_t, n_t) \in C([0,T]; L^q(\mathbb{R}^3) \cap L^2(\mathbb{R}^3))
,\\
&u\in C([0,T];D^2(\mathbb{R}^3)\cap D_0^1(\mathbb{R}^3))
\cap L^2(0,T;D^{2,q}(\mathbb{R}^3))
,\\
&u_t\in L^2\big(0,T;D_0^1(\mathbb{R}^3)\big),\qquad
\sqrt{\rho}\,u_t\in L^\infty\big(0,T;L^2(\mathbb{R}^3)\big),
\end{align*}
for some $q\in(3,6]$, and $(\rho,n,u)$ satisfies system \eqref{twofluid}
almost everywhere in $\mathbb{R}^3\times(0,T)$ together with the initial conditions (\ref{non-initial}) and (\ref{non-cauchy}).
\end{definition}

Now we are in the position to state our main result in the paper.
\begin{theorem}\label{theorem 1.2}
Let the initial data $(\rho_0, n_0, u_0)$ satisfy
\[
0\leq\rho_0\le\max\limits_{\mathbb{R}^3}\rho_0:=\bar\rho,\quad 0\leq n_0\le\max\limits_{\mathbb{R}^3}n_0:=\bar n,
\]
\[
\rho_0,n_0\in W^{1,q}(\mathbb{R}^3)\cap H^1(\mathbb{R}^3),\quad
u_0\in D^2(\mathbb{R}^3)\cap D_0^1(\mathbb{R}^3),
\]
for some positive constants $\bar\rho , \,\bar n$, and some $q\in(3,6]$, and that the initial total energy is finite, namely,
\[
E_0=\int\Big(\frac12\rho_0|u_0|^2
+\frac{1}{\gamma-1}\rho_0^\gamma
+\frac{1}{\Gamma-1}n_0^\Gamma\Big) < \infty,\quad \gamma>1,\quad \Gamma>1,
\]
and  the following compatibility condition
\beq\label{compatibility}
-\mu\Delta u_0-(\mu+\lambda)\nabla\di  u_0+\nabla P(\rho_0,n_0)=\sqrt{\rho_0}g,
\eeq holds for some $g\in L^2(\mathbb{R}^3)$.
Then there exists a positive constant $\varepsilon_0,$  depending only on $\mu,\lambda,\gamma$ and  $\Gamma$
 such that the initial-value problem (\ref{twofluid}) and (\ref{pressure2}), (\ref{non-initial}) and (\ref{non-cauchy}) admits a unique global strong solution provided that
\be\label{isen-smallness}
\bar P^{\frac{ 3}{\gamma}} E_0 \big(1+\bar P^{\frac{ 3}{\gamma}}\bar PE_0^2\big)\big(\|\nabla u_0\|_{L^2}^2+\|P_0\|_{L^2}^2\big)\le \varepsilon_0,
\ee where   $P_0=P(\rho_0,n_0)$ and $\bar P=\bar\rho^\gamma+\bar n^\Gamma$.
\end{theorem}

\begin{remark}
The quantity on the left-hand side of (\ref{isen-smallness}) is scaling invariant for (\ref{scl-inva-non}). 
\end{remark}

\section{Proof of Theorem \ref{theorem 1.2}}\label{sec 2}

We first show the local existence theory  of a unique strong solution with vacuum  to the Cauchy problem for system \eqref{twofluid} and (\ref{pressure2}), which can be established by slightly modifying the proof in \cite{Wen-yaozhu} due to different pressure law.
\begin{lemma}[Local existence]\label{local}
Assume the initial data $(\rho_0, n_0, u_0)$ satisfy Theorem \ref{theorem 1.2}.
Then there exists a small time $T_0>0$ and a unique strong solution $(\rho,n,u)$ to the Cauchy problem
\eqref{twofluid} and (\ref{pressure2}), \eqref{non-initial}, \eqref{non-cauchy} such that
\begin{align*}
&
(\rho,\, n) \in C\big([0,T_0]; W^{1,q}(\mathbb{R}^3)\cap H^1(\mathbb{R}^3)\big), \\
&(\rho_t,\, n_t) \in C\big([0,T_0]; L^q(\mathbb{R}^3)\cap L^2(\mathbb{R}^3)\big), \\
&u \in C\big([0,T_0]; D^2(\mathbb{R}^3)\cap D^1_0(\mathbb{R}^3)\big)
      \cap L^2\big(0,T_0; D^{2,q}(\mathbb{R}^3)\big), \\
&u_t \in L^2\big(0,T_0; D^1_0(\mathbb{R}^3)\big), \
\sqrt{\rho}\, u_t \in L^\infty\big(0,T_0; L^2(\mathbb{R}^3)\big),
\end{align*}
where $q\in(3,6]$.
\end{lemma}

The global existence of the solution will be obtained by using the standard continuity argument together with a series of global-in-time a priori estimates. Next, we state the following Gagliardo-Nirenberg inequality that will be used frequently throughout the rest of the paper.
\begin{lemma}[\cite{GN}]\label{GN}
Let $1 \leq r,p,q \leq \infty$, $m \in \mathbb{N}^{+}$, $k \in\mathbb{N}$, with $0 \leq k < m$, and let $\theta, r$ be such that
\begin{equation}
0 \leq \theta \leq 1 - k/m
\end{equation}
and
\begin{equation}
(1-\theta)\left(\frac{1}{p} - \frac{m - k}{3}\right) + \theta\left(\frac{1}{q} + \frac{k}{3}\right) = \frac{1}{r} .
\end{equation}
Then there exists a constant $C = C(m, p, q, \theta, k) > 0$ such that
\begin{equation}\label{2.1}
\|\nabla^k f\|_{L^r(\mathbb{R}^3)} \leq C\|f\|_{L^q(\mathbb{R}^3)}^\theta \|\nabla^m f\|_{L^p(\mathbb{R}^3)}^{1-\theta}
\end{equation}
for every $f \in L^q(\mathbb{R}^3) \cap D^{m,p}(\mathbb{R}^3)$.
\end{lemma}

The following Zlotnik inequality will be applied to derive the uniform-in-time upper bound of the density $\rho$ and $n$.
\begin{lemma}[\cite{Z2000}]\label{zolonik}
Let the function y satisfies
\bex
y^\prime(t)\leq g(y)+b^\prime(t)\,\,\mathrm{on}\,[0,T],\, y(0)=y^0,
\eex
with $g\in C(R)$ and $y,b\in W^{1,1}(0,T)$. If $g(\infty)=-\infty$ and
\begin{equation}
b(t_2)-b(t_1)\le N_0+N_1(t_2-t_1)
\end{equation}
for all $0\le t_1<t_2\le T$ with some $N_0\ge0$ and $N_1\ge0$, then
\begin{equation}\label{2.4}
y(t)\le \max\{y^0,\bar{\zeta}\}+N_0<\infty\, \mathrm{on}\, [0,T],
\end{equation} where $\bar{\zeta}$ is a constant such that
\bex
g(\zeta)\le -N_1\, \mathrm{for}\, \zeta\ge\bar{\zeta}.
\eex
\end{lemma}

\subsection{Global-in-time {\it a priori} estimates}

\begin{proposition}\label{prop1}
Under the hypotheses of Theorem~\ref{theorem 1.2}, letting
$(\rho,\,n,\,u)$ be a strong solution to \eqref{twofluid} and (\ref{pressure2}), \eqref{non-initial}, \eqref{non-cauchy}
on $\mathbb{R}^3\times[0,T]$ for $T>0$, and defining the energy functional
\begin{equation}\label{energy}
\mathcal E(t)
:=\bar P^{\frac{ 3}{\gamma}} E_0
\Big(1+\bar P^{\frac{3}{\gamma}}\bar P E_0^2\Big)
\Big(\|\nabla u(t)\|_{L^2}^2+\|P(t)\|_{L^2}^2\Big),
\end{equation}
there exists a small positive constant $\varepsilon$ independent of $T$,
$\bar\rho$, $\bar n$ and the initial data and determined by \eqref{min240}, \eqref{smallness2}, and \eqref{smallness3}, such that if
\begin{equation}\label{aprio-assum}
\begin{cases}
P(t)\le 4{\bar P},\quad\mathrm{for}\,\, \mathrm{any}\,\,t\in[0,T],\\[2mm]
\displaystyle
\sup_{t\in[0,T]}\mathcal E(t)\le 2\varepsilon,
\end{cases}
\end{equation}
then the following improved estimates
\begin{equation}
\begin{cases}
P(t)\le 3{\bar P}, \quad\mathrm{for}\,\, \mathrm{any}\,\,t\in[0,T],\\[2mm]
\displaystyle
\sup_{t\in[0,T]}\mathcal E(t)\le \dfrac32\,\varepsilon,
\end{cases}
\end{equation} hold, provided that the initial data satisfy
\bex
\bar{P}^{\frac{ 3}{\gamma}} E_0 \big(1+\bar P^{\frac{ 3}{\gamma}}\bar PE_0^2\big)\big(\|\nabla u_0\|_{L^2}^2+\|P_0\|_{L^2}^2\big)\le \varepsilon_0
\eex where $\varepsilon_0=\frac{3\varepsilon}{2C}$ for a positive constant $C$ determined by Corollary \ref{isen-cor:2.6}.
\end{proposition}

The proof of Proposition \ref{prop1} consists of Lemma \ref{isen-le:2.8} and Corollary \ref{isen-cor2}. Throughout the rest of the paper, let $C>1$ denote a generic positive constant independent of the initial data, $T$, $\varepsilon_0$, $\varepsilon$, $\bar\rho$ and $\bar n$.
\begin{lemma}\label{le2.1}Under the hypotheses of Proposition \ref{prop1}, there holds
\begin{equation}\label{result1}
\begin{split}
&\int \Big(\frac{1}{2}\rho|u|^2+\frac{1}{\gamma-1}\rho{^\gamma}+\frac{1}{\Gamma-1}n{^\Gamma}\Big)+\int_0^t\int\big[\mu|\nabla u|^2+(\mu+\lambda)|\di u|^2\big]\,ds\\&= \int \Big(\frac{1}{2}\rho_0|u_0|^2+\frac{1}{\gamma-1}\rho_0{^\gamma}+\frac{1}{\Gamma-1}n_0{^\Gamma}\Big)\doteq E_0
\end{split}
\end{equation} for any $t\in[0,T]$.
\end{lemma}
\pf
The energy equality is standard, since the solution is smooth enough. More specifically, multiplying \eqref{twofluid}$_3$ by $u$,  integrating by parts over $\mathbb{R}^3$, and using the continuity equation, \eqref{result1} can be obtained.
\endpf

\begin{lemma}\label{isen-le: 2.3}Under the hypotheses of Proposition \ref{prop1}, there holds
\be\label{i-lenau}
\begin{split}
&\int|\nabla u|^2+\int_0^t\int \rho|\dot u|^2\,ds\\ \le& C\int P^{2}+C\int(|\nabla u_0|^2+P_0^{2})+
C{\bar P^{\frac{3}{\gamma}}} \int_0^t\|\nabla u\|_{L^2}^4\Big(\|\nabla u\|_{L^2}^2+\|P\|_{L^2}^2\Big)\,ds\\&+C{\bar P^{\frac{1}{\gamma}}} \int_0^t\|\nabla u\|_{L^2}^2\|P\|_{L^3}^2\,ds
\end{split}
\ee
 for any $t\in[0,T]$.
\end{lemma}
\pf
Multiplying the momentum equation \eqref{twofluid}$_3$ by $u_t$, and integrating  over $\mathbb{R}^3$, we obtain
\be\label{211}
\begin{split}
&\int \rho |\dot u|^2+\frac{1}{2}\dt \int \left(\mu|\nabla u|^2+(\mu+\lambda)|\di u|^2\right)\\=&\int \rho u\cdot\nabla u \cdot \dot u+\dt \int P\di u-\int P_t\di u.
\end{split}
\ee
Introduce the effective viscous flux
\bex
G=(2\mu+\lambda)\di  u-P.
\eex
Using  $\di  u=\frac{G+P}{2\mu+\lambda}$, we rewrite
\[
-\int P_t\di  u
=
-\frac{1}{2(2\mu+\lambda)}\dt \int P^2
-\frac{1}{2\mu+\lambda}\int P_t\,G .
\]
Substituting the above equation  into  \eqref{211} yields
\beq\label{212}
\begin{split}
&\int \rho|\dot u|^2
+\frac12\dt \int\left(\mu|\nabla\times u|^2+\frac{1}{2\mu+\lambda}|G|^2\right)\\
=\;&
\int \rho u\cdot\nabla u\cdot\dot u
-\frac{1}{2\mu+\lambda}\int P_t\,G .
\end{split}
\eeq
By using the continuity equations for $\rho$ and $n$, we start by multiplying the first equation in \eqref{twofluid} by \( \Gamma n^{\Gamma - 1} \) and the second equation by \( \gamma \rho^{\gamma - 1} \), which leads to the following
\beq\label{2.14}
\begin{split}
P_t+u\cdot \nabla P+\gamma P\di u+(\Gamma-\gamma)n^\Gamma\di u=0.
\end{split}
\eeq
For simplicity, we write it as follows
\beq\label{213}
\begin{split}
P_t+\di (Pu)=\mathcal{R},
\end{split}
\eeq
where $\mathcal{R}=-(\gamma-1)P\di u-(\Gamma-\gamma)n^\Gamma \di u$.  Substituting \eqref{213} into \eqref{212} and integrating by parts yields the decomposition
\[
-\int P_t\,G
=
-\int u\, P \cdot\nabla G
-\int \mathcal{R}\,G .
\]
Consequently, using the estimates above, we obtain
\beq
\begin{split}
&\int \rho|\dot u|^2
+\frac12\dt \int\left(\mu|\nabla\times u|^2+\frac{1}{2\mu+\lambda}|G|^2\right)\\
\le\;&\frac12\int \rho|\dot u|^2+C\int \rho|u|^2|\nabla u|^2
+C\int |u|| P||\nabla G|+C\int |\di u|| P|| G|+C\int |\di u|| n^\Gamma|| G|.\\
\end{split}
\eeq
Then, by absorbing \(\frac12 \int \rho|\dot u|^2\) by the left-hand side, we further have
\beq
\begin{split}\label{2.5}
&\int \rho|\dot u|^2
+\dt \int\left(\mu|\nabla\times u|^2+\frac{1}{2\mu+\lambda}|G|^2\right)\\
\le\;&C\int \rho|u|^2|\nabla u|^2
+C\int |u|| P||\nabla G|+C\int |\di u|| P|| G|+C\int |\di u|| n^\Gamma|| G|\\
\le\;& C\int \rho|u|^2|\nabla u|^2
+C\int |u|| P||\nabla G|+C\int |\nabla u|| P|| G|\\
\nonumber=:\;& I_1+I_2+I_3.
\end{split}
\eeq
 By Sobolev inequality and the boundedness of $\rho$, we have
\beq\label{I1}
\begin{split}
I_1&\le C\int \rho |u|^2|\nabla u|^2
\\ &\le C\|\rho\|_{L^\infty} \|u\|_{L^6}^2\|\nabla u\|_{L^3}^2
\\ &\le C\bar P^{\frac{1}{\gamma}}\|\nabla u\|_{L^2}^2\Big(\|\nabla\times u\|_{L^3}^2+\|G\|_{L^3}^2+\|P\|_{L^3}^2\Big).
\end{split}
\eeq
Taking $\mathrm{curl}$ and $\di $ on both sides of (\ref{twofluid})$_3$, respectively, we have
\be\label{i-eqGw}
\begin{cases}
\mu\Delta(\mathrm{curl}u)=\mathrm{curl}(\rho \dot u),\\[2mm]
\Delta G=\di (\rho \dot u).
\end{cases}
\ee
Thus, by standard elliptic estimates
\beq\label{na Gw}
\begin{split}\|\nabla G\|_{L^{2}}+\|\nabla\nabla\times u\|_{L^2} \le
C{\bar P^{\frac{1}{2\gamma}}}\|\sqrt{\rho} \dot u\|_{L^2}.
\end{split}
\eeq
Substituting (\ref{na Gw}) into (\ref{I1}), we have
\beq\label{I1+1}
\begin{split}
I_1&\le C\bar P^\frac{3}{2\gamma} \|\nabla u\|_{L^2}^2\Big(\|\nabla u\|_{L^2}+\|P\|_{L^2}\Big)\|\sqrt{\rho}\dot u\|_{L^2}+C\bar P^{\frac{1}{\gamma}} \|\nabla u\|_{L^2}^2\|P\|_{L^3}^2.
\end{split}
\eeq
By using the H\"older inequality, Sobolev inequality, and (\ref{na Gw}), we have
\beq\label{i-I 2}\begin{split}
I_2+I_3&\le C\int | P| |u||\nabla G|+C\int |\nabla u|| P|| G|
\\ &\le C\|u\|_{L^6}\|P\|_{L^3}\|\nabla G\|_{L^2}+ C\|G\|_{L^6}\|P\|_{L^3}\|\nabla u\|_{L^2}
\\ &\le C\|\nabla u\|_{L^2}\|P\|_{L^3}{\bar P^{\frac{1}{2\gamma}}}\|\sqrt{\rho} \dot u\|_{L^2}.
\end{split}
\eeq
Combining the estimates for $I_1,\,I_2$ and $I_3,$  and applying Young's inequality, we obtain
\beq\begin{split}\label{2.11}
&\int \rho |\dot u|^2+\dt \int\left(\mu|\nabla \times u|^2+\frac{1}{2\mu+\lambda}|G|^2\right)
\\ &\le C {\bar P^{\frac{3}{\gamma}}} \|\nabla u\|_{L^2}^4\Big(\|\nabla u\|_{L^2}^2+\|P\|_{L^2}^2\Big)+C{\bar P^{\frac{1}{\gamma}}} \|\nabla u\|_{L^2}^2\|P\|_{L^3}^2.
\end{split}
\eeq
Integrating (\ref{2.11}) over $(0,t)$ yields the desired estimate  (\ref{i-lenau}).

\endpf

\begin{lemma}\label{isen-le2.4}
Under the hypotheses of Proposition \ref{prop1}, there holds
\begin{align}\label{le:1}
\|(-\Delta)^{-1}\di (\rho u)\|_{L^\infty}  \le C\bar P^\frac{3}{4\gamma}E_0^\frac{1}{4}\|\nabla u\|_{L^2}^\frac{1}{2}
\end{align} for any $t\in[0,T]$.
\end{lemma}
\pf As in \cite{Wen2025}, Using  Gagliardo-Nirenberg inequality, Sobolev inequality  and (\ref{result1}), we have
\begin{align}\label{i-Lsecond}
\nonumber\|(-\Delta)^{-1}\di (\rho u)\|_{L^\infty}\le& C\|(-\Delta)^{-1}\di (\rho u)\|_{L^6}^\frac{1}{2}\|\nabla(-\Delta)^{-1}\di (\rho u)\|_{L^6}^\frac{1}{2}
\\ \le& \nonumber C\|\rho u\|_{L^2}^\frac{1}{2}\|\rho u\|_{L^6}^\frac{1}{2}
\\ \le&\nonumber C\|\rho\|^\frac{3}{4}_{L^\infty}\|\sqrt{\rho} u\|_{L^2}^\frac{1}{2}\|\nabla u\|_{L^2}^\frac{1}{2}
\\ \le&\nonumber C \bar P^\frac{3}{4\gamma}\|\sqrt{\rho} u\|_{L^2}^\frac{1}{2}\|\nabla u\|_{L^2}^\frac{1}{2}
\\ \le&  C\bar P^\frac{3}{4\gamma}E_0^\frac{1}{4}\|\nabla u\|_{L^2}^\frac{1}{2}.
\end{align}
\endpf

\begin{lemma}\label{isen-le:2.5}Under the hypotheses of Proposition \ref{prop1}, there holds
\begin{align}\label{i-dtrhop+11}
 \nonumber\int P^{2}+\int_0^t\int P^{3}\,ds &\le C\int P_0^{2}+ C\bar P^{\frac{ 3}{\gamma}}E_0\sup\limits_{s\in(0,t)}\big(\|P\|_{L^2}^2+\|\nabla u(s)\|_{L^2}^2\big) \\
&\times \Big(\int_0^t\|\sqrt{\rho}\dot u\|_{L^2}^2\,ds+\sup\limits_{s\in(0,t)}\|\nabla u(s)\|_{L^2}^2\Big)
\end{align} for any $t\in[0,T]$.
\end{lemma}
\pf
Without loss of generality, we assume that $\Gamma\ge\gamma$.  Multiplying (\ref{2.14}) by $2P$,  we have
\beq\label{2.21}
\begin{split}
(P^{2})_t+\nabla\cdot (uP^2)+(2\gamma-1)P^2\di u+2(\Gamma-\gamma)Pn^\Gamma \di u=0.
\end{split}
\eeq
Integrating the above identity over $\mathbb{R}^3$ yields
\beq\label{rho infty 11}
\begin{split}
\dt\int P^{2}
=-(2\gamma-1)\int P^{2}\di u-2(\Gamma-\gamma)\int n^\Gamma P\di u.
\end{split}
\eeq
Using  the definition of  $\di  u$, we derive
\beq\label{rho infty 11}
\begin{split}
&\dt \int P^{2}+\frac{(2\gamma-1)}{2\mu+\lambda}\int P^{3}+\frac{2(\Gamma-\gamma)}{2\mu+\lambda}\int n^\Gamma P^2\\
&=\frac{-(2\gamma-1)}{2\mu+\lambda}\int P^{2}G+\frac{-2(\Gamma-\gamma)}{2\mu+\lambda}\int n^\Gamma PG.
\end{split}
\eeq
Substituting the expression of $G$ given in $(\ref{i-eqGw})$, namely,
 \be\label{i-equation of G+11}
G=-(-\Delta)^{-1}\di (\rho u)_t-(-\Delta)^{-1}\di \di (\rho u\otimes u)
\ee
into the right-hand side of \eqref{rho infty 11}, we obtain
\beq
\begin{split}\label{2.22}
&\dt \int P^2 + \frac{2\gamma - 1}{2\mu + \lambda} \int P^3 +\frac{2(\Gamma-\gamma)}{2\mu+\lambda}\int n^\Gamma P^2\\
&-\frac{2\gamma-1}{2\mu+\lambda}\dt \int P^{2}(-\Delta)^{-1}\di (\rho u)-\frac{2(\Gamma-\gamma)}{2\mu+\lambda}\dt \int n^\Gamma P(-\Delta)^{-1}\di (\rho u)
 \\=& -\frac{2\gamma-1}{2\mu+\lambda}\int (P^{2})_t(-\Delta)^{-1}\di (\rho u)
+\frac{2\gamma-1}{2\mu+\lambda}\int P^{2}(-\Delta)^{-1}\di \di (\rho u\otimes u)\\ &
-\frac{2(\Gamma-\gamma)}{2\mu+\lambda}\int (n^\Gamma P)_t(-\Delta)^{-1}\di (\rho u)
+\frac{2(\Gamma-\gamma)}{2\mu+\lambda}\int n^\Gamma P (-\Delta)^{-1}\di \di (\rho u\otimes u)
\\=&:J_1+J_2+J_3+J_4.
\end{split}
\eeq
We now estimate the terms $J_1$--$J_4$ on the right-hand side of \eqref{2.22}.
For $J_1,$ by  using the H\"{o}lder inequality and \eqref{2.21} , we have
\beq\label{2.26}
\begin{split}
J_1&=-\frac{2\gamma-1}{2\mu+\lambda}\int (P^{2})_t(-\Delta)^{-1}\di (\rho u)\\
&\leq C\|P\|_{L^3}^2\|u\|_{L^6}\|\rho u\|_{L^6}+C\|P\|_{L^3}^2\|\di u\|_{L^3}\|(-\Delta)^{-1}\di (\rho u)\|_{L^\infty}.
\end{split}
\eeq
For $J_2$ and $J_4,$  applying  the H\"{o}lder and Sobolev embedding inequality, we obtain
\beq\label{2.27}
\begin{split}
J_2&=\frac{2\gamma-1}{2\mu+\lambda}\int P^{2}(-\Delta)^{-1}\di \di (\rho u\otimes u)\\
&\leq C\|P\|_{L^3}^2\|(-\Delta)^{-1}\di \di (\rho u\otimes u)\|_{L^3}\\
&\leq C\|P\|_{L^3}^2\|\rho u\otimes u\|_{L^3}\\
&\leq C\| \rho\|_{L^\infty}\|P\|_{L^3}^2\|u\|_{L^6}^2\\
&\leq C\bar P^{\frac{1}{\gamma}}\|P\|_{L^3}^2\|\nabla u\|_{L^2}^2.
\end{split}
\eeq
Following a similar argument as in  \eqref{2.27}, we can estimate \(J_4\) as
\beq\label{2.27+}
\begin{split}
J_4&=\frac{2(\Gamma-\gamma)}{2\mu+\lambda}\int n^\Gamma P (-\Delta)^{-1}\di \di (\rho u\otimes u)\\
&\leq C\bar P^{\frac{1}{\gamma}}\|P\|_{L^3}^2\|\nabla u\|_{L^2}^2.
\end{split}
\eeq
For $J_3,\, $we first multiply \eqref{twofluid}$_1$ by $\Gamma n^{\Gamma-1}$ and  obtain
\beq\label{n_t}
\begin{split}
(n^\Gamma)_t+\nabla\cdot (un^\Gamma)+(\Gamma-1)n^\Gamma\di u=0.
\end{split}
\eeq
Using \eqref{2.14} and \eqref{n_t}, we have
\beq
\begin{split}
(n^\Gamma P)_t&=(n^\Gamma)_t P+n^\Gamma P_t\\
&=[-\nabla\cdot (un^\Gamma)-(\Gamma-1)n^\Gamma\di u]P+n^\Gamma[-u\cdot \nabla P-\gamma P\di u-(\Gamma-\gamma)n^\Gamma\di u]\\
&=-u\cdot\nabla n^\Gamma P-\Gamma n^\Gamma\di u P- u\cdot\nabla P n^\Gamma -\gamma n^\Gamma\di uP-(\Gamma-\gamma)(n^\Gamma)^2\di u\\
&=-u\cdot\nabla (Pn^\Gamma)-(\Gamma+\gamma)n^\Gamma\di u P-(\Gamma-\gamma)(n^\Gamma)^2\di u\\
&=-\nabla\cdot (Pn^\Gamma u)-(\Gamma+\gamma-1)n^\Gamma\di u P-(\Gamma-\gamma)(n^\Gamma)^2\di u.
\end{split}
\eeq
Therefore, using the H\"{o}lder inequality and similar arguments as for $J_1$, we have
\beq\label{2.30}
\begin{split}
J_3=&-\frac{2(\Gamma-\gamma)}{2\mu+\lambda}\int (n^\Gamma P)_t(-\Delta)^{-1}\di (\rho u)\\
\leq& C\|P\|_{L^3}\|n^\Gamma\|_{L^3}\|u\|_{L^6}\|\rho u\|_{L^6}+C\|P\|_{L^3}\|n^\Gamma\|_{L^3}\|\di u\|_{L^3}\|(-\Delta)^{-1}\di (\rho u)\|_{L^\infty}\\
&+C\|n^\Gamma\|_{L^3}^2\|\di u\|_{L^3}\|(-\Delta)^{-1}\di (\rho u)\|_{L^\infty}\\
\leq& C\|P\|_{L^3}^2\|u\|_{L^6}\|\rho u\|_{L^6}+C\|P\|_{L^3}^2\|\di u\|_{L^3}\|(-\Delta)^{-1}\di (\rho u)\|_{L^\infty}.
\end{split}
\eeq
Substituting \eqref{2.26}, \eqref{2.27}, \eqref{2.27+} and \eqref{2.30} to (\ref{2.22}), we obtain
\be\label{rho2ga}
\begin{split}
&\dt \int P^2 + \frac{2\gamma - 1}{2\mu + \lambda} \int P^3 +\frac{2(\Gamma-\gamma)}{2\mu+\lambda}\int n^\Gamma P^2\\
&-\frac{2\gamma-1}{2\mu+\lambda}\dt \int P^{2}(-\Delta)^{-1}\di (\rho u)-\frac{2(\Gamma-\gamma)}{2\mu+\lambda}\dt \int n^\Gamma P(-\Delta)^{-1}\di (\rho u)
\\ \le&C\|P\|_{L^3}^2\|u\|_{L^6}\|\rho u\|_{L^6}+C\|P\|_{L^3}^2\|\di u\|_{L^3}\|(-\Delta)^{-1}\di (\rho u)\|_{L^\infty}\\&+C\bar P^{\frac{1}{\gamma}}\|P\|_{L^3}^2\|\nabla u\|_{L^2}^2.
\end{split}
\ee
By virtue of the definition of  $\di u$ and \eqref{na Gw}, we get
\be\label{divuL3}
\begin{split}
\|\di u\|_{L^3}\le&C\|G\|_{L^2}^\frac{1}{2}\|G\|_{L^6}^\frac{1}{2}+C\|P\|_{L^3}
\\ \le&C\big(\|\nabla u\|_{L^2}^\frac{1}{2}+\|P\|_{L^2}^\frac{1}{2}\big)\|\nabla G\|_{L^2}^\frac{1}{2}+C\|P\|_{L^3}
\\ \le&C{\bar P^{\frac{1}{4\gamma}}}\big(\|\nabla u\|_{L^2}^\frac{1}{2}+\|P\|_{L^2}^\frac{1}{2}\big)\|\sqrt{\rho}\dot u\|_{L^2}^\frac{1}{2}+C\|P\|_{L^3},
\end{split}
\ee
where we have used Sobolev inequality, the H\"{o}lder inequality,  and the definition of the effective viscous flux $G$.
Putting (\ref{divuL3}) to (\ref{rho2ga}), and using (\ref{le:1}), we have
\be\label{rho2ga+1}
\begin{split}
&\dt \int P^2 + \frac{2\gamma - 1}{2\mu + \lambda} \int P^3 +\frac{2(\Gamma-\gamma)}{2\mu+\lambda}\int n^\Gamma P^2\\
&-\frac{2\gamma-1}{2\mu+\lambda}\dt \int P^{2}(-\Delta)^{-1}\di (\rho u)-\frac{2(\Gamma-\gamma)}{2\mu+\lambda}\dt \int n^\Gamma P(-\Delta)^{-1}\di (\rho u)
\\ \le&C\bar P^{\frac{1}{\gamma}}\|P\|_{L^3}^2\|\nabla u\|_{L^2}^2+C\bar P^{\frac{1}{\gamma}}\|P\|_{L^3}^2\|\sqrt{\rho}\dot u\|_{L^2}^\frac{1}{2}\|\sqrt{\rho}u\|_{L^2}^\frac{1}{2}\|\nabla u\|_{L^2}\\&+C\bar P^{\frac{1}{\gamma}}\|P\|_{L^3}^2\|P\|_{L^2}^\frac{1}{2}\|\sqrt{\rho}\dot u\|_{L^2}^\frac{1}{2}\|\sqrt{\rho}u\|_{L^2}^\frac{1}{2}
\|\nabla u\|_{L^2}^\frac{1}{2}\\&+C\|P\|_{L^3}^3\bar P^\frac{3}{4\gamma}\|\sqrt{\rho}u\|_{L^2}^\frac{1}{2}\|\nabla u\|_{L^2}^\frac{1}{2}
\\ \le&C\bar P^{\frac{1}{\gamma}}\|P\|_{L^3}^2\|\nabla u\|_{L^2}^2+C\bar P^{\frac{1}{\gamma}}\|P\|_{L^3}^2\|\sqrt{\rho}\dot u\|_{L^2}^\frac{1}{2}E_0^\frac{1}{4}\|\nabla u\|_{L^2}\\&+C\bar P^{\frac{1}{\gamma}}\|P\|_{L^3}^2\|P\|_{L^2}^\frac{1}{2}\|\sqrt{\rho}\dot u\|_{L^2}^\frac{1}{2}E_0^\frac{1}{4}\|\nabla u\|_{L^2}^\frac{1}{2}\\&+C\|P\|_{L^3}^3\bar P^\frac{3}{4\gamma}E_0^\frac{1}{4}\|\nabla u\|_{L^2}^\frac{1}{2}.
\end{split}
\ee
Integrating (\ref{rho2ga+1}) over $(0,t)$, and using the Cauchy inequality, (\ref{aprio-assum}), and (\ref{le:1}), we have
\be\label{rho2ga+2}
\begin{split}
&\int P^{2}+\frac{2\gamma-1}{2(2\mu+\lambda)}\int_0^t\int P^{3}\,ds+\frac{2(\Gamma-\gamma)}{2(2\mu+\lambda)}\int_0^t\int n^\Gamma P^2\,ds
\\ \le&C\int P_0^{2}+C\bar P^{\frac{ 3}{\gamma}}\int_0^t\|\nabla u\|_{L^2}^6\,ds+C\bar P^{\frac{ 3}{\gamma}}\int_0^t\|\sqrt{\rho}\dot u\|_{L^2}^\frac{3}{2}E_0^\frac{3}{4}\|\nabla u\|_{L^2}^3\,ds\\&+C\bar P^{\frac{ 3}{\gamma}}\int_0^t\|P\|_{L^2}^\frac{3}{2}\|\sqrt{\rho}\dot u\|_{L^2}^\frac{3}{2}E_0^\frac{3}{4}\|\nabla u\|_{L^2}^\frac{3}{2}\,ds\\
\le&C\int P_0^{2}+K_1+K_2+K_3,
\end{split}
\ee
where we have used the smallness assumption \eqref{aprio-assum} such that
\beq\label{min240}
\begin{split}
C\bar P^\frac{3}{4\gamma}\|\sqrt{\rho}u\|_{L^2}^\frac{1}{2}\|\nabla u\|_{L^2}^\frac{1}{2}\leq C(2\epsilon)^\frac{1}{4}\leq \min\{\frac{2\gamma-1}{2(2\mu+\lambda)},\frac{2\mu+\lambda}{2(2\gamma-1 )}\}.
\end{split}
\eeq
For $K_i$ (i=1,2,3), using (\ref{result1}), we have
\be\label{II1-410}
\begin{split}
K_1 \le C\bar P^{\frac{ 3}{\gamma}}E_0\sup\limits_{s\in(0,t)}\big(\|\nabla u(s)\|_{L^2}^2\big)\sup\limits_{s\in(0,t)}\|\nabla u(s)\|_{L^2}^2,
\end{split}
\ee
\be\label{II1-41}
\begin{split}
K_2\le& C\bar P^{\frac{ 3}{\gamma}}E_0^\frac{3}{4}\sup\limits_{s\in(0,t)}\big(\|\nabla u(s)\|_{L^2}^2\big)\Big(\int_0^t\|\sqrt{\rho}\dot u\|_{L^2}^2\,ds\Big)^\frac{3}{4}\Big(\int_0^t\|\nabla u\|_{L^2}^4\,ds\Big)^\frac{1}{4}\\
\le& C\bar P^{\frac{ 3}{\gamma}}E_0^\frac{3}{4}\sup\limits_{s\in(0,t)}\big(\|\nabla u(s)\|_{L^2}^2\big)
\sup\limits_{s\in(0,t)}\|\nabla u(s)\|_{L^2}^\frac{1}{2}\Big(\int_0^t\|\sqrt{\rho}\dot u\|_{L^2}^2\,ds\Big)^\frac{3}{4}\Big(\int_0^t\|\nabla u\|_{L^2}^2\,ds\Big)^\frac{1}{4}\\
\le& C\bar P^{\frac{ 3}{\gamma}} E_0\sup\limits_{s\in(0,t)}\big(\|\nabla u(s)\|_{L^2}^2\big)\Big(\int_0^t\|\sqrt{\rho}\dot u\|_{L^2}^2\,ds+\sup\limits_{s\in(0,t)}\|\nabla u(s)\|_{L^2}^2\Big),
\end{split}
\ee
and
\be\label{II1-4}
\begin{split}
K_3\le&C\bar P^{\frac{ 3}{\gamma}}E_0^\frac{3}{4}\sup\limits_{s\in(0,t)}\big(\|P\|_{L^2}^\frac{3}{2}\|\nabla u(s)\|_{L^2}^\frac{1}{2}\big)\Big(\int_0^t\|\sqrt{\rho}\dot u\|_{L^2}^2\,ds\Big)^\frac{3}{4}\Big(\int_0^t\|\nabla u\|_{L^2}^4\,ds\Big)^\frac{1}{4}
\\ \le&C\bar P^{\frac{ 3}{\gamma}}E_0\sup\limits_{s\in(0,t)}\big(\|P\|_{L^2}^\frac{3}{2}\|\nabla u(s)\|_{L^2}^\frac{1}{2}\big)\Big(\int_0^t\|\sqrt{\rho}\dot u\|_{L^2}^2\,ds+\sup\limits_{s\in(0,t)}\|\nabla u(s)\|_{L^2}^2\Big).
\end{split}
\ee
Putting (\ref{II1-410}), (\ref{II1-41}) and (\ref{II1-4}) to (\ref{rho2ga+2}), we get
\begin{align*}
 \nonumber&\int P^{2}+\int_0^t\int P^{3}\,ds\\ \le &C\int P_0^{2}+ C\bar P^{\frac{ 3}{\gamma}}E_0\sup\limits_{s\in(0,t)}\big(\|P\|_{L^2}^2+\|\nabla u(s)\|_{L^2}^2\big)\Big(\int_0^t\|\sqrt{\rho}\dot u\|_{L^2}^2\,ds+\sup\limits_{s\in(0,t)}\|\nabla u(s)\|_{L^2}^2\Big).
\end{align*}
The proof of Lemma \ref{isen-le:2.5} is complete.
\endpf

\begin{corollary}\label{isen-cor1}
Under the hypotheses of Proposition \ref{prop1}, there holds
\be\label{cor1}
\int|\nabla u|^2+\int_0^t\int \rho |\dot u|^2\,ds \le C\int\left(|\nabla u_0|^2+P_0^{2}\right) +\frac{1}{4}\int_0^t\|P\|_{L^3}^3\,ds
\ee for any $t\in[0,T]$.
\end{corollary}
\pf
Using (\ref{result1}), (\ref{i-lenau}), and the H\"older inequality, we have
\be\label{nauL2}
\begin{split}
&\int|\nabla u|^2+\int_0^t\int  \rho |\dot u|^2\,ds\\ \le& C\int P^{2}+C\int\left(|\nabla u_0|^2+P_0^{2}\right)+C\bar P^{\frac{ 3}{\gamma}} E_0\sup\limits_{s\in(0,t)}\Big[\|\nabla u(s)\|_{L^2}^2\big(\|\nabla u(s)\|_{L^2}^2+\|P\|_{L^2}^2\big)\Big]\\&+C\bar P^{\frac{1}{\gamma}} \sup\limits_{s\in(0,t)}\|\nabla u(s)\|_{L^2}^\frac{4}{3}E_0^\frac{1}{3}\Big(\int_0^t\|P\|_{L^3}^3\,ds\Big)^\frac{2}{3}.
\end{split}
\ee
Next, using (\ref{i-dtrhop+11}) to control $\displaystyle\int P^2$ in (\ref{nauL2}), we obtain
\be\label{nauL2+1}
\begin{split}
&\sup\limits_{s\in(0,t)}\int|\nabla u(s)|^2+\int_0^t\int \rho |\dot u|^2\,ds\\
\le& C\int\left(|\nabla u_0|^2+P_0^{2}\right)+C\bar P^{\frac{ 1}{\gamma}} \sup\limits_{s\in(0,t)}\|\nabla u(s)\|_{L^2}^\frac{4}{3}E_0^\frac{1}{3}\Big(\int_0^t\|P\|_{L^3}^3\,ds\Big)^\frac{2}{3}\\
\le& C\int\left(|\nabla u_0|^2+P_0^{2}\right)+C\bar P^{\frac{ 3}{\gamma}} E_0 \sup\limits_{s\in(0,t)}\|\nabla u(s)\|_{L^2}^2\sup\limits_{s\in(0,t)}\|\nabla u(s)\|_{L^2}^2 +\frac{1}{8}\int_0^t\|P\|_{L^3}^3\,ds\\
\le& C\int\left(|\nabla u_0|^2+P_0^{2}\right)+\frac{1}{2}\sup\limits_{s\in(0,t)}\|\nabla u(s)\|_{L^2}^2 +\frac{1}{8}\int_0^t\|P\|_{L^3}^3\,ds,
\end{split}
\ee
where we have used Young's inequality together with the a priori assumption (\ref{aprio-assum})
\begin{equation}\label{smallness2}
C\bar P^{\frac{ 3}{\gamma}} E_0
\sup_{s\in(0,t)} \big( \|P(s)\|_{L^2}^2 + \|\nabla u(s)\|_{L^2}^2 \big)
\;\le\; 2C \varepsilon \;\le\; \frac{1}{4},
\end{equation}
which ensures that the nonlinear terms can be absorbed by the left-hand side of \eqref{nauL2+1}.
\endpf

\begin{corollary}\label{isen-cor:2.6}
Under the hypotheses of Proposition \ref{prop1}, there holds
\be\label{cor2}\begin{split}
\int(|\nabla u|^2+P^{2})+\int_0^t\int \rho |\dot u|^2\,ds+\int_0^t\int P^{3}\,ds
 \le C\int\left(|\nabla u_0|^2+P_0^{2}\right)
\end{split}
\ee for any $t\in[0,T]$.
\end{corollary}
\pf We use the energy functional \(\mathcal{E}(t)\) defined in \eqref{energy}
$$
\mathcal E(t)
:=\bar P^{\frac{ 3}{\gamma}} E_0
\Big(1+\bar P^{\frac{ 3}{\gamma}}\bar P E_0^2\Big)
\Big(\|\nabla u(t)\|_{L^2}^2+\|P(t)\|_{L^2}^2\Big).
$$
Putting (\ref{i-dtrhop+11}) and (\ref{cor1}) together and using (\ref{smallness2}), we have
\begin{equation}\label{cor3}
\begin{split}
&\int (|\nabla u|^2+P^{2})
+\int_0^t\int\rho |\dot u|^2\,ds
+\int_0^t\int P^{3}\,ds  \\
\le\;&
C\int\left(|\nabla u_0|^2+P_0^{2}\right)
+\frac{1}{4}\int_0^t\|P\|_{L^3}^3\,ds\\
&+ C\mathcal E(t)\Big(\int_0^t\|\sqrt{\rho}\dot u\|_{L^2}^2\,ds+\sup\limits_{s\in(0,t)}\|\nabla u(s)\|_{L^2}^2\Big).
\end{split}
\end{equation}
The last two terms on the right-hand side of \eqref{cor3} can be absorbed by the left-hand side. Hence, \eqref{cor2} follows immediately.
\endpf

\begin{lemma}\label{isen-le:2.8}
Under the hypotheses of Proposition \ref{prop1}, there holds
\bex
P(x,t)\le 3\bar{P}
\eex for any $(x,t)\in\mathbb{R}^3\times[0,T]$.
\end{lemma}
\pf
Define the particle trajectories $X(s;x,t)$  given by
    $$ \left\{
      \begin{array}{l}
       \displaystyle\frac{d}{ds}X(s;x,t)=u\left(X(s;x,t), s\right),\ 0\le s<t,\\[2mm]
        X(t; x,t)=x.\\
      \end{array}
      \right.
    $$
    As in \cite{Wen-ZhuSIMA}, letting $(x,t) \in \mathbb{R}^3 \times [0,T]$ be fixed, and
denoting
$$
\rho^\delta(y,s)=\rho(y,s)+\delta\exp\{-\int_0^s\di u\left(X(\tau;x,t),\tau\right)\,d\tau\}>0,
$$
it is easy to verify that
$$
\frac{d}{d s}\rho^\delta\left(X(s;x,t),
s\right)+\rho^\delta\left(X(s;x,t),
s\right)\di u\left(X(s;x,t),s\right)=0.
$$
This yields
\be\label{i-rhoidentity}\begin{split} Y_1^\prime(s)\leq g_1(s)+b^\prime(s),
\end{split}\ee where
\bex
\begin{split}
&Y_1(s)=\ln\rho^\delta\left(X(s;x,t),s\right),\\
&g_1(s)=-\frac{\rho^\gamma\left(X(s;x,t),s\right)}{2\mu+\lambda},\\
&b(s)=-\frac{1}{2\mu+\lambda}\int_0^sG\left(X(\tau;x,t),\tau\right)\,d\tau,
\end{split}
\eex
and $G=(2\mu+\lambda)\di u-P$.
For $G$, we have
\bex\begin{split} G\left(X(t;x,\tau),\tau\right)=&-\frac{d}{d\tau}[(-\Delta)^{-1} \di (\rho
u)]+[u_i,R_{ij}](\rho u_j),
\end{split}\eex
where $[u_i,R_{ij}]=u_iR_{ij}-R_{ij}u_i$ and
$R_{ij}=\partial_i(-\Delta)^{-1}\partial_j$. This deduces
\bex
\begin{split}
b(t_2)-b(t_1)=&\frac{1}{2\mu+\lambda}(-\Delta)^{-1}
\di (\rho u)(t_2)-\frac{1}{2\mu+\lambda}(-\Delta)^{-1}
\di (\rho
u)(t_1)\\&-\frac{1}{2\mu+\lambda}\int_{t_1}^{t_2}[u_i,R_{ij}](\rho u_j)\,d\tau\\
\le& \frac{2}{2\mu+\lambda}\sup\limits_{0\le t\le T}\|(-\Delta)^{-1} \di (\rho
u)(t)\|_{L^\infty}+\frac{1}{2\mu+\lambda}\int_{t_1}^{t_2}\|[u_i,R_{ij}](\rho
u_j)\|_{L^\infty}\,d\tau\\
=&II_1+II_2.
\end{split}
\eex
For $II_1$, it follows from Lemma \ref{isen-le2.4} that
\be\label{II_1}
\begin{split}
II_1\le  C\sup\limits_{0\le t\le T}\bar P^\frac{3}{4\gamma}E_0^\frac{1}{4}\|\nabla u\|_{L^2}^\frac{1}{2}.
\end{split}
\ee
For $II_2$, it follows from  $(33)$ in \cite{Desjardin, Wen-ZhuSIMA} and \eqref{na Gw} that
\be \label{II2}
\begin{split}
II_2\le& C\int_{t_1}^{t_2}\bar P^\frac{1}{\gamma}\|\nabla u\|_{L^2}\|\nabla u\|_{L^6}\,d\tau
\\ \le&C\int_{t_1}^{t_2}\bar P^\frac{1}{\gamma}\|\nabla u\|_{L^2}\Big(\bar P^\frac{1}{2\gamma}\|\sqrt{\rho}\dot u\|_{L^2}+\|P\|_{L^6}\Big)\,d\tau
\\ \le&C\bar P^\frac{1}{\gamma}\Big(\int_{t_1}^{t_2}\|\nabla u\|_{L^2}^2\,d\tau\Big)^\frac{1}{2}\Big(\int_{t_1}^{t_2}\bar P^\frac{1}{\gamma}\|\sqrt{\rho}\dot u\|_{L^2}^2\,d\tau\Big)^\frac{1}{2}\\&
+C\bar P^\frac{1}{\gamma}\Big(\int_{t_1}^{t_2}\|\nabla u\|_{L^2}^2\,d\tau\Big)^\frac{1}{2}\Big(\int_{t_1}^{t_2}\|P\|_{L^6}^2\,d\tau\Big)^\frac{1}{2}.
\end{split}
\ee
Note that the following estimate holds
\be\label{III-0}
\begin{split}
&\bar P^\frac{1}{\gamma}\Big(\int_{t_1}^{t_2}\|\nabla u\|_{L^2}^2\,d\tau\Big)^\frac{1}{2}\Big(\int_{t_1}^{t_2}\|P\|_{L^6}^2\,d\tau\Big)^\frac{1}{2}\\
\le&
\bar P^\frac{1}{\gamma}\Big(\int_{t_1}^{t_2}\|\nabla u\|_{L^2}^2\,d\tau\Big)^\frac{1}{2}\bar P^\frac{1}{2}\Big(\int_{t_1}^{t_2}\|P\|_{L^3}\,d\tau\Big)^\frac{1}{2}
\\ \le&\bar P^\frac{1}{\gamma}\Big(\int_{t_1}^{t_2}\|\nabla u\|_{L^2}^2\,d\tau\Big)^\frac{1}{2}\bar P^\frac{1}{2}\Big(\int_{t_1}^{t_2}\|P\|_{L^3}^3\,d\tau\Big)^\frac{1}{6}(t_2-t_1)^\frac{1}{3}
\\ \le&\bar P^\frac{1}{\gamma}\Big(\int_{t_1}^{t_2}\|\nabla u\|_{L^2}^2\,d\tau\Big)^\frac{1}{2}\bar P^\frac{1}{2}\Big[\int\left(|\nabla u_0|^2+P_0^{2}\right)\Big]^\frac{1}{6}(t_2-t_1)^\frac{1}{3}
\\ \le&C\bar P^\frac{3}{2\gamma}\Big(\int_{t_1}^{t_2}\|\nabla u\|_{L^2}^2\,d\tau\Big)^\frac{3}{4}\bar P^\frac{1}{4}\Big[\int\left(|\nabla u_0|^2+P_0^{2}\right)\Big]^\frac{1}{4}+\frac{\bar P}{2\mu+\lambda}(t_2-t_1),
\end{split}
\ee
where we have used the H\"older inequality, Young's inequality and (\ref{cor2}).

As in \cite{Wen2025}, we have
\bex
\begin{split}
b(t_2)-b(t_1)\leq&N_0+N_1(t_2-t_1),
\end{split}
\eex
where
\begin{equation*}
\begin{split}
N_0=&C\sup\limits_{0\le t\le T}\bar P^\frac{3}{4\gamma}E_0^\frac{1}{4}\|\nabla u\|_{L^2}^\frac{1}{2}\\
&+C\bar P^\frac{3}{2\gamma} E_0^\frac{1}{2}\Big[\int\left(|\nabla u_0|^2+P_0^{2}\right)\Big]^\frac{1}{2}
+C\bar P^\frac{3}{2\gamma}E_0^\frac{3}{4}\bar P^\frac{1}{4}\Big[\int\left(|\nabla u_0|^2+P_0^{2}\right)\Big]^\frac{1}{4},\\[2mm]
N_1=&\frac{\bar P}{2\mu+\lambda}.
\end{split}
\end{equation*}Note that the new quantity $\bar P$ in $N_0$ is very essential to guarantee that the smallness of the initial data is scaling invariant.}

Recalling that
$$
g_1=g_1(Y_1):=-\frac{\Big(e^{Y_1}
-\delta\exp\{\displaystyle-\int_0^s\di u\left(X(\tau;x,t),\tau\right)\,d\tau\}\Big)^\gamma}{2\mu+\lambda},
$$
it is clear that $g_1(\infty)=-\infty$. For a fixed $\delta>0$, define
\begin{equation}
\bar{Y}_1^\delta:=\ln\Big\{\bar P^\frac{1}{\gamma}+\delta\exp\{\int_0^T\|\di u(\cdot,\tau)\|_{L^\infty}\,d\tau\}\Big\},
\end{equation}
and letting $Y_1\ge\bar{Y}_1^\delta$, there holds that
\begin{equation}
\begin{split}
g_1(Y_1)=&-\frac{\Big(e^{Y_1}
-\delta\exp\{\displaystyle-\int_0^s\di u\left(X(\tau;x,t),\tau\right)\,d\tau\}\Big)^\gamma}{2\mu+\lambda}
\\ \le&-\frac{\Big(e^{Y_1}
-\delta\exp\{\displaystyle\int_0^T\|\di u(\cdot,\tau)\|_{L^\infty}\,d\tau\}\Big)^\gamma}{2\mu+\lambda}
\\ \le&-N_1.
\end{split}
\end{equation}
Applying Lemma \ref{zolonik} along the trajectory, it follows that
\begin{equation*}
\begin{split}
&\rho^\delta\left(x,s\right)\le \max\Big\{\bar\rho+\delta,\,\exp\{\bar{Y}_1^\delta\}\Big\}\exp\{N_0\}.
\end{split}
\end{equation*}
Letting $\delta$ tend to zero, we have
\begin{equation}\label{rho upp}
\begin{split}
&\rho\left(x,s\right)\le \bar P^\frac{1}{\gamma}\exp\{N_0\}.
\end{split}
\end{equation}
Similarly, for $n$,  it holds that
$$
n^\delta(y,s)=n(y,s)+\delta\exp\{-\int_0^s\di u\left(X(\tau;x,t),\tau\right)\,d\tau\}>0.
$$
It is easy to verify that
$$
\frac{d}{d s}n^\delta\left(X(s;x,t),
s\right)+n^\delta\left(X(s;x,t),
s\right)\di u\left(X(s;x,t),s\right)=0.
$$
We can obtain the result by the same method
\bex
\begin{split}
&Y_2(s)=\ln n^\delta\left(X(s;x,t),s\right),\\
&g_2(s)=-\frac{n^\Gamma\left(X(s;x,t),s\right)}{2\mu+\lambda},\\
&b(s)=-\frac{1}{2\mu+\lambda}\int_0^sG\left(X(\tau;x,t),\tau\right)\,d\tau.
\end{split}
\eex

Similarly, we have
$$
g_2=g_2(Y_2):=-\frac{\Big(e^{Y_2}
-\delta\exp\{\displaystyle-\int_0^s\di u\left(X(\tau;x,t),\tau\right)\,d\tau\}\Big)^\Gamma}{2\mu+\lambda},
$$
it is clear that $g_2(\infty)=-\infty$. For a fixed $\delta>0$, define
\begin{equation}
\bar{Y}_2^\delta:=\ln\Big\{\bar P^\frac{1}{\Gamma}+\delta\exp\{\int_0^T\|\di u(\cdot,\tau)\|_{L^\infty}\,d\tau\}\Big\},
\end{equation}
and letting $Y_2\ge\bar{Y}_2^\delta$, there holds that
\begin{equation}
\begin{split}
g_2(Y_2)=&-\frac{\Big(e^{Y_2}
-\delta\exp\{\displaystyle-\int_0^s\di u\left(X(\tau;x,t),\tau\right)\,d\tau\}\Big)^\Gamma}{2\mu+\lambda}
\\ \le&-\frac{\Big(e^{Y_2}
-\delta\exp\{\displaystyle\int_0^T\|\di u(\cdot,\tau)\|_{L^\infty}\,d\tau\}\Big)^\Gamma}{2\mu+\lambda}
\\ \le&-N_1.
\end{split}
\end{equation}
By virtue of  Lemma \ref{zolonik} applied along the trajectory, we derive that
\begin{equation*}
\begin{split}
n^\delta\left(x,s\right)\le \max\Big\{\bar n+\delta,\,\exp\{\bar{Y}_2^\delta\}\Big\}\exp\{N_0\}.
\end{split}
\end{equation*}
Letting $\delta$ tend to zero, we have
\begin{equation}\label{n upp}
\begin{split}
n\left(x,s\right)\le \bar P^\frac{1}{\Gamma}\exp\{N_0\}.
\end{split}
\end{equation}
It is obvious  that
\begin{equation}\label{P upp}
\begin{split}
P\left(x,s\right)\le \bar P\Big(\exp\{\gamma N_0\}+\exp\{\Gamma N_0\}\Big).
\end{split}
\end{equation}
Letting
\be\label{smallness3}
\begin{split}
&(\gamma+\Gamma) N_0\\
\le&C\sup\limits_{0\le t\le T}\bar P^\frac{3}{4\gamma}E_0^\frac{1}{4}\|\nabla u\|_{L^2}^\frac{1}{2}\\
&+C\bar P^\frac{3}{2\gamma} E_0^\frac{1}{2}\Big[\int\left(|\nabla u_0|^2+P_0^{2}\right)\Big]^\frac{1}{2}
+C\bar P^\frac{3}{2\gamma}E_0^\frac{3}{4}\bar P^\frac{1}{4}\Big[\int\left(|\nabla u_0|^2+P_0^{2}\right)\Big]^\frac{1}{4}
\\ \le& C(2\varepsilon)^\frac{1}{4}+C(2\varepsilon)^\frac{1}{2}\le \ln\frac{3}{2},
\end{split}
\ee
and using (\ref{P upp}), we finish the proof Lemma \ref{isen-le:2.8}.
\endpf

\begin{corollary}\label{isen-cor2}
Under the hypotheses of Proposition \ref{prop1}, there holds
\begin{equation}
\mathcal E(t)\le \frac{3}{2}\,\varepsilon
\end{equation}
for any $t\in[0,T]$.
\end{corollary}

\pf
By Corollary \ref{isen-cor:2.6} we have the uniform-in-time estimate
\bex
\|\nabla u(t)\|_{L^2}^2 + \|P(t)\|_{L^2}^2
\le C \big(\|\nabla u_0\|_{L^2}^2 + \|P_0\|_{L^2}^2\big),
\eex for any $t\in[0,T]$.
This together with the definition of the  energy functional \(\mathcal{E}(t)\) in \eqref{energy} yields
\[
\mathcal{E}(t)
\le C\bar P^\frac{3}{\gamma} E_0 \Big(1+\bar P^\frac{3}{\gamma} \bar P E_0^2\Big) \big(\|\nabla u_0\|_{L^2}^2 + \|P_0\|_{L^2}^2\big)
= C \mathcal{E}(0),
\]
where \(C>0\) is the constant determined by Corollary \ref{isen-cor:2.6}.
To guarantee \(\mathcal{E}(t)\le \frac32 \varepsilon\), it suffices to choose \(\varepsilon_0\) satisfying
\begin{equation}\label{def:epsilon0}
\mathcal{E}(0) \le \varepsilon_0 := \frac{3\varepsilon}{2C},
\end{equation}
we have
\[
\mathcal{E}(t) \le C \mathcal{E}(0)  \le \frac{3}{2}\,\varepsilon.
\]

This finishes the proof of Corollary \ref{isen-cor2}.
\endpf

\subsection{Proof of Theorem \ref{theorem 1.2}}

Let $(\rho,n,u)$ be the local strong solution stated in Lemma \ref{local} on $[0,T^*)$ where $T^*$ is the maximal existence time for the strong solution. By virtue of Proposition \ref{prop1}, Corollary \ref{isen-cor2}, and the continuity method, the following estimates
\begin{equation}\label{concl}
\begin{cases}
P(t)\le 4{\bar P},\quad\mathrm{for}\,\, \mathrm{any}\,\,t\in [0,T^*),\\[2mm]
\displaystyle
\sup_{t\in[0,T^*)}\mathcal E(t)\le 2\varepsilon,
\end{cases}
\end{equation}
hold, provided that the initial energy satisfies
\[
\mathcal{E}(0) \le \varepsilon_0,
\]
where $\varepsilon_0$ is given by (\ref{def:epsilon0}).
In particular, (\ref{cor2}), Sobolev inequality, and (\ref{concl})$_1$ yield
\be\label{uL6}
\|\rho\|_{L^\infty(0,T^*; L^\infty)}+\|n\|_{L^\infty(0,T^*; L^\infty)}+\|\sqrt{\rho}u\|_{L^\infty(0,T^*; L^6)}\le C.
\ee
Since $T^*$ is the maximal existence time, Remark 1.3 in \cite{Wen-yaozhu} demonstrates that if $T^*<\infty$ and the pressure function in (\ref{twofluid}) is replaced by (\ref{pressure1}), the following blowup criterion
\be\label{blowup criterion}
\limsup\limits_{T\rightarrow T^*}\Big(\|\rho\|_{L^\infty(0,T; L^\infty)}+\|\sqrt{\rho}u\|_{L^\infty(0,T; L^6)}\Big)=\infty
\ee holds with initial densities satisfying
\be\label{density condition}
0\le \underline{s}_0\rho_0\le n_0\le \overline{s}_0\rho_0
\ee for some constants $\underline{s}_0>0$ and $\overline{s}_0>0$.
In fact, if $\rho$ in (\ref{blowup criterion}) is replaced by
$(\rho,n)$ and the pressure function is considered in this paper, the domination condition of density (\ref{density condition}) can be removed with slightly modification. Therefore it contradicts with (\ref{uL6}) if $T^*<\infty$. There must be that $T^*=\infty$. The proof of Theorem \ref{theorem 1.2} is complete.

\endpf

\section*{Acknowledgements}
This work was supported by the National Natural Science Foundation of China $\#$ 12471209 and by Guangzhou Basic and Applied
Research Projects \#SL2024A04J01206.


\begin{thebibliography}{00}
\bibitem{BarrettLuSuli2017} J.W. Barrett, Y. Lu and E. S\"{u}li, {\em Existence of large-data finite-energy global weak solutions to a compressible Oldroyd-B model,} Commun. Math. Sci., 15 (2017), 1265--1323.
\bibitem{Desjardin}
B. Desjardins, {\em Regularity of weak solutions of the compressible
isentropic Navier-Stokes equations,} Commun. Partial Differential Equations, 22(1997), 977--1008.

\bibitem{Evje2017} S. Evje, {\em An integrative multiphase model for cancer cell migration under influence of physical cues from the tumor microenvironment,} Chem. Eng. Sci., 165 (2017), 240--259.

\bibitem{EK2008} S. Evje, K.H.  Karlsen, {\em Global existence of weak solutions for a viscous two-phase model,} J. Differ. Equations, 245 (2008), 2660--2703.

\bibitem{EvjeWen2018} S. Evje, H.Y. Wen, {\em A Stokes two-fluid model for cell migration that can account for physical cues in the microenvironment,} SIAM J. Math. Anal., 50 (2018), 86--118.

\bibitem{EvjeWenZhu2017} S. Evje, H.Y. Wen and C.J. Zhu, {\em On global solutions to the viscous liquid-gas model with unconstrained transition to single-phase flow,} Math. Models Methods Appl. Sci., 27 (2017), 323--346.

\bibitem{GuoYangYao2011} Z.H. Guo, J. Yang and L. Yao, {\em Global strong solution for a three-dimensional viscous liquid-gas two-phase flow model with vacuum,} J. Math. Phys., 52 (2011), 093102.



\bibitem{HaoLi2012} C.C. Hao, H.L. Li, {\em Well-posedness for a multidimensional viscous liquid-gas two-phase flow model,} SIAM J. Math. Anal., 44 (2012), 1304--1332.

\bibitem{Hu} X.P. Hu, {\em Global existence for two dimensional compressible magnetohydrodynamic flows with zero magnetic diffusivity,} arXiv:1405.0274v1.
\bibitem{jiangzhang2017}
S. Jiang, J.W. Zhang, {\em On the non-resistive limit and the magnetic boundary-layer for one-dimensional compressible magnetohydrodynamics,} Nonlinearity, 30 (2017), 3587--3612.


\bibitem{KMN2024} M. Kalousek, S. Mitra, $\mathrm{\breve{S}}$. Ne$\mathrm{\breve{c}}$asov\'a, {\em The existence of a weak solution for a compressible
multicomponent fluid structure interaction problem,} J. Math. Pures Appl., 184 (2024), 118--189.

\bibitem{GN}
G. Leoni, {\em A first course in Sobolev spaces}, second edition,
Graduate Studies in Mathematics, 181, Amer. Math. Soc., Providence, RI, 2017.



\bibitem{LiSun2019} Y. Li, Y.Z. Sun, {\em Global weak solutions to a two-dimensional compressible MHD equations of viscous non-resistive fluids,} J. Differ. Equations, 267 (2019), 3827--3851.





\bibitem{NP2020} A. Novotn\'y, M. Pokorn\'y, {\em Weak solutions for some compressible multicomponent fluid models,} Arch. Ration. Mech. Anal., 235 (2020), 355--403.
\bibitem{Oldroyd} J.G. Oldroyd, {\em Non-Newtonian effects in steady motion of some idealized elastico-viscous liquids,} Proc. Roy. Soc. London Ser. A, 245 (1958), 278--297.

\bibitem{VWY-2019} A.F. Vasseur, H.Y. Wen and C. Yu, {\em Global weak solution to the viscous two-fluid model with finite energy,} J. Math. Pures Appl., 125 (2019), 247--282.


\bibitem{W2021} H.Y. Wen, {\em On global solutions to a viscous compressible two-fluid model with unconstrained transition to single-phase flow in three dimensions,} Calc. Var. Partial Differential Equations, 60 (2021), Paper No. 158, 38 pp.

\bibitem{Wen2025} H.Y. Wen, {\em Global well-posedness of compressible Navier-Stokes equations with vacuum and smallness on scaling invariant quantity in $\mathbb{R}^3$,} Adv. Math., 482 (2025), 110628.
\bibitem{Wen-yaozhu} H.Y. Wen, L. Yao and C.J. Zhu, {\em A blow-up criterion of strong solution to a 3D viscous liquid-gas two-phase flow model with vacuum,} J. Math. Pures Appl., 97 (2012), no.~3, 204--229.

\bibitem{Wen-ZhuSIMA} H.Y. Wen, C.J. Zhu, {\em Global solutions to the three-dimensional full compressible Navier-Stokes equations with vacuum at infinity in some classes of large data,} SIAM J. Math. Anal., 49 (2017), 162--221.
\bibitem{Wen-Zhucpaa}
H.Y. Wen,  C.J.  Zhu, \emph{ Remarks on global weak solutions to a two-fluid type model,} Commun. Pure Appl. Anal.,  20 (2021), no.~7-8, 2839--2856.


\bibitem{Wuwu2017} J.H. Wu, Y.F. Wu, {\em Global small solutions to the compressible 2D magnetohydrodynamic system without magnetic diffusion,} Adv. Math., 310 (2017), 759--888.

\bibitem{Wuzhu2022} J.H. Wu, Y. Zhu, {\em Global well-posedness for 2D non-resistive compressible MHD system in periodic domain,} J. Funct. Anal., 283 (2022), Paper No. 109602, 49.

\bibitem{Yao-Zhu2}
L. Yao, C.J. Zhu, {\em Existence and uniqueness of global weak solution
to a two-phase flow model with vacuum,} Math. Ann., 349 (2011), 903--928.

\bibitem{Yu2021} H.B. Yu, {\em Global strong solutions to the 3D viscous liquid-gas two-phase flow model,} J. Differ. Equations. 272 (2021), 732--759.




\bibitem{Z2000} A.A. Zlotnik, {\em Uniform estimates and stabilization of symmetric solutions of a system of quasilinear equations,} Diff. Uravn., 36 (2000), 634--646; translation in Differ. Equ., 36 (2000), 701--716.
\end{thebibliography}
\end{document}